\documentclass[a4paper,12pt,final]{article}
\pdfoutput=1
\usepackage{geometry}
\usepackage[T1]{fontenc}
\usepackage{ifthen}
\usepackage{amsmath}
\usepackage{amsfonts}
\usepackage{amsthm}
\usepackage{amssymb}
\usepackage{enumerate}
\usepackage[all]{xy}
\usepackage{rotating}
\usepackage{epsfig}
\usepackage{color}
\usepackage{hyperref}
\hypersetup{colorlinks=true,linkcolor=black,citecolor=black}
\newtheorem{Lem}{Lemma}[section]
\newtheorem{Con}[Lem]{Construction}
\newtheorem{Pro}[Lem]{Proposition}
\newtheorem{Cor}[Lem]{Corollary}
\newtheorem{Thm}[Lem]{Theorem}
\newtheorem*{Thm*}{Theorem}
\theoremstyle{definition}
\newtheorem{Exa}[Lem]{Example}
\newtheorem*{Exa*}{Example}
\newtheorem{Def}[Lem]{Definition}
\newtheorem{Rem}[Lem]{Remark}
\newtheorem*{Rem*}{Remark}
\newcommand{\s}[1]{\left(\begin{smallmatrix}#1\end{smallmatrix}\right)}
\begin{document}
\thispagestyle{empty}
\begin{center} 
\textbf{\Large Duality of non-exposed faces}\\
\vspace{.3cm}
Stephan Weis\footnote{\texttt{sweis@mis.mpg.de}}\\
Max Planck Institute for Mathematics in the Sciences\\
Leipzig, Germany\\
\vspace{.1cm}
November 27, 2011
\end{center}
\noindent
{\small\textbf{\emph{Abstract -- }}
Given any polar pair of convex bodies we study its conjugate face maps
and we characterize conjugate faces of non-exposed faces in terms of
normal cones. The analysis is carried out using the positive hull
operator which defines lattice isomorphisms linking three Galois
connections. One of them assigns conjugate faces between the convex
bodies. The second and third Galois connection is defined between
the touching cones and the faces of each convex body separately.
While the former is well-known, we introduce the latter in this
article for any convex set in any finite dimension. We demonstrate
our results about conjugate faces with planar convex bodies and planar 
self-dual convex bodies, for which we also include constructions.}\\[5mm]
%
%
%
%
%
%
\section{Non-exposed faces and dual convex bodies}
\par
Duality of faces of a dual pair of closed convex cones was studied in 
\cite{Barker_Dual} with regard to the lattice of the inclusion ordering.
This duality corresponds to the conjugate face map between faces of a
polar pair of convex bodies. In this article we study the restriction
of the conjugate face map to non-exposed faces. E.g.\ it will become
clear that a face which is conjugate to a non-exposed face is singular
(its normal cone has at least dimension two). We prove that such faces
are fully characterized by a so-called {\it incomplete} normal cone. 
\par
Incomplete normal cones of planar convex bodies have a simple description
by so-called {\it mixed} and {\it free corners}. Examples are given in 
Figure~\ref{fig:dual_pics} and~\ref{fig:self_pics}. The conjugate face
map restricts to a surjective map from the non-exposed points of a
planar convex body onto the mixed and free corners of the polar convex
body
\begin{equation}
\label{eq:easy_map}\textstyle
\xymatrix{
\{\text{ {\it non-exposed points} }\}\quad \ar@{->>}[rr]
 && \quad\{\text{ mixed corners and free corners }\}\,.
}
\end{equation}
\par
The idea underlying this article is to use (\ref{eq:easy_map}), 
and its generalization in any dimension, to study non-exposed faces
of a projection of the state space of the matrix algebra
${\rm Mat}(N,\mathbb C)$. The polar convex body of a projection 
is an affine section of that state space, see \S2.4 in 
\cite{Weis_support}. Its singular points (with incomplete normal cone)
may be studied by analyzing an associated determinantal variety, using 
techniques of algebraic geometry. Our interest in non-exposed points
of projected state spaces lies in quantum information theory, they
seem to cause discontinuities in certain information measures 
\cite{Knauf_Weis}.
\par
Planar projections of state spaces are studied in operator theory 
under the name of numerical range, see e.g.\ \cite{Dunkl} and the 
references therein. The question when numerical range has non-exposed 
points was solved in \cite{Rodman} for $N=3$. Recently, numerical
range was studied in \cite{Henrion} from the point of view of convex
algebraic geometry whose aim is to use techniques from algebraic
geometry for studying convex semialgebraic sets. Important examples
of such sets are spectrahedra which generalize the state space of 
${\rm Mat}(N,\mathbb C)$ and which are popular in optimization.
Current questions in the field are concerned with convex duality and 
non-exposed faces, see e.g.\ \cite{Netzer,Sanyal}. Our interest in 
self-dual convex bodies is influenced by the present discussion of 
self-duality in the axiomatic foundations of quantum theory 
\cite{Janotta,Mueller,Wilce}.
\par
This article is organized as follows. Constructions for dual convex 
bodies and a general construction for planar self-dual convex bodies
are explained in \S\ref{sec:examples}. A Galois connection between
touching cones and faces of a convex set is defined in \S\ref{sec:galois}.
In \S\ref{sec:pairs_of_bodies} we study conjugate faces of any polar
pair of convex bodies. We demonstrate our results in \S\ref{sec:2D}
with planar convex bodies. In particular we give a general construction 
for planar self-dual convex bodies without non-exposed points.
%
%
%
%
\section{Constructions of dual convex bodies}
\label{sec:examples}
\par
We introduce constructions for dual and self-dual convex bodies 
(mainly in dimension two). They are used to generate examples to
demonstrate non-exposed points and their relation to the singular
points studied in \S\ref{sec:2D}.
\par
In the $n$-dimensional Euclidean vector space
$(\mathbb R^n,\langle\cdot,\cdot\rangle)$ we denote the norm of 
$u\in\mathbb R^n$ by $|u|:=\sqrt{\langle u,u\rangle}$. In $\mathbb R^n$ we 
shall use the standard scalar product. The {\it polar} of a subset 
$C\subset\mathbb R^n$ is 
$C^\circ:=\{u\in\mathbb R^n\mid\langle u,v\rangle\leq 1\;\forall v\in C\}$ 
and the {\it dual} of $C$ is $C^*:=\{u\in\mathbb R^n\mid 
1+\langle u,v\rangle\geq 0\;\forall v\in C\}=-C^\circ$. The subset
$C\subset\mathbb R^n$ is {\it self-dual} if $C^*=C$. We denote the
{\it interior} of $C$ by ${\rm int}(C)$ and its boundary by
$\partial(C):=C\setminus{\rm int}(C)$. 
\par
The first construction is Corollary~16.5.2 in \cite{Rockafellar}:
\begin{Con}
For any family $\{C_i\}_{i\in I}$ of convex sets in $\mathbb R^n$ ($I$
is an index set) we have
\begin{align}
\label{eq:1}\textstyle
(\text{convex hull of }\{C_i\mid i\in I\})^*&
\;=\;
\bigcap\{C_i^*\mid i\in I\}\,.
\end{align} 
\end{Con}
\begin{Exa}
\label{ex:1}
The convex set $C\subset\mathbb R^2$ 
depicted in Figure~\ref{fig:dual_pics} c) is the convex hull of the unit
disk $D:=\{u\in\mathbb R^2\mid |u|\leq 1\}$ and of the point $\s{0\\2}$. 
We have $D^*=D$ (e.g.\ using (\ref{eq:1})) and
$\{\s{0\\2}\}^*=\{(x,y)\in\mathbb R^2\mid y\geq-\tfrac 1{2}\}$. The dual 
$C^*=D^*\cap\{\s{0\\2}\}^*$ is depicted in Figure~\ref{fig:dual_pics} d).
\end{Exa}
\par
In the sequel let $K\subset\mathbb R^n$ denote a {\it convex body}, i.e.\
a convex and compact subset, and let $0\in{\rm int}(K)$. 
By Theorem~1.6.1 in \cite{Schneider} the polar $K^\circ$ is a convex body
with $0\in{\rm int}(K^\circ)$ and $(K^\circ)^\circ=K$. Obviously the dual
$K^*$ is a convex body with $0\in{\rm int}(K^*)$ and $(K^*)^*=K$.
A second construction for the dual convex body arises from the 
{\it support function} of a convex $C\subset\mathbb R^n$ in the
direction $u\in\mathbb R^n$,
\[\textstyle
h_C(u)
\;:=\;
\sup\{\langle x,u\rangle\mid x\in C\}\,.
\]
The {\it radial function} of the convex body $K$ is
\[\textstyle
\rho_K(u)
\;:=\;
\sup\{\lambda\geq 0\mid\lambda u\in K\}\,.
\]
Theorem~1.7.6 in \cite{Schneider} shows for all $u\in\mathbb R^n$
that $\rho_{K^\circ}(u)=1/h_K(u)$ holds, hence
\begin{equation}
\label{eq:rad_supp}\textstyle
\rho_{K^*}(u)
\;=\;
1/h_K(-u)\,.
\end{equation}
This equation includes $\rho_{K^*}(0)=\infty$ and $h_K(0)=0$ with the
convention of $1/0=\infty$. 
\begin{Con}
\label{con:support_radial}
The boundary of the dual convex body ${K^*}$ is parametrized from the
unit sphere by the support function of $K$,
\[\textstyle
S^{n-1}:=\{u\in\mathbb R^n\mid|u|=1\}\;\to\;\partial K^*\,,
\quad
u\;\mapsto\;\rho_{K^*}(u)u=u/h_K(-u)\,.
\]
\end{Con}
{\em Proof:\/}
The map $u\mapsto\rho_{K^*}(u)u$ defined on the unit sphere $S^{n-1}$ 
extends to a positively homogeneous function $\mathbb R^n\to\mathbb R^n$
by setting $0\mapsto 0$ and $u\mapsto\rho_{K^*}(\tfrac u{|u|})u$
for $u\neq 0$. The Theorem of Sz.~Nagy (see e.g.\ \S{}VIII.1 in 
\cite{Berge}) 
shows that this function, called {\it radial projection}, is a 
homeomorphism between the unit ball and ${K^*}$. In particular,
$S^{n-1}\to\partial{K^*}$, $u\mapsto\rho_{K^*}(u)u$ is a
parametrization of the boundary of ${K^*}$. The radial function of $K^*$
is expressed by the support function of $K$ in (\ref{eq:rad_supp}).
\hspace*{\fill}$\Box$\\
\begin{figure}
\centerline{%
\begin{picture}(2.14,3.2)
\put(0,0){\includegraphics[height=3.2cm, bb=0 0 267 400, 
clip=]{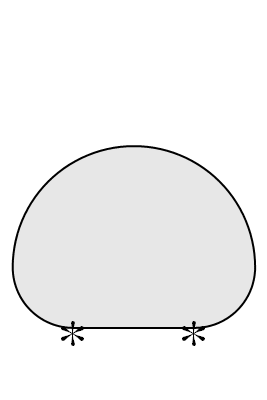}}
\put(0,0){a)}
\end{picture}
\hspace{1.5cm}
\begin{picture}(2.14,3.2)
\put(0,0){\includegraphics[height=3.2cm, bb=0 0 267 400, 
clip=]{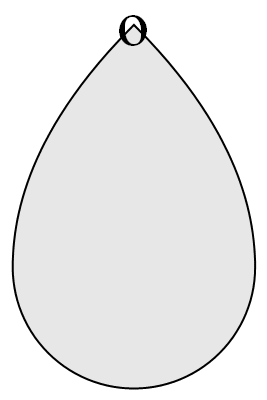}}
\put(0,0){b)}
\end{picture}
\hspace{1.5cm}
\begin{picture}(2.14,3.2)
\put(0,0){\includegraphics[height=3.2cm, bb=0 0 267 400, 
clip=]{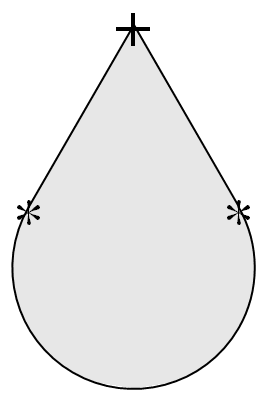}}
\put(0,0){c)}
\end{picture}
\hspace{1.5cm}
\begin{picture}(2.14,3.2)
\put(0,0){\includegraphics[height=3.2cm, bb=0 0 267 400, 
clip=]{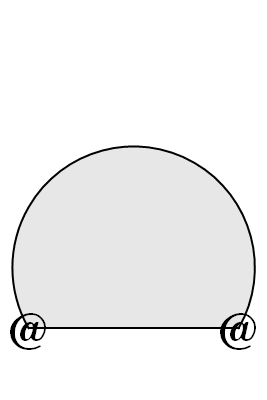}}
\put(0,0){d)}
\end{picture}}
\caption{\label{fig:dual_pics} 
The convex sets a) and b) likewise c) and d) are duals of each other. 
Markings indicate non-exposed points (*), polyhedral corners 
(+), mixed corners (@) and free corners (o). The drawings 
have equal scaling and their origin is aligned vertically.}
\end{figure}
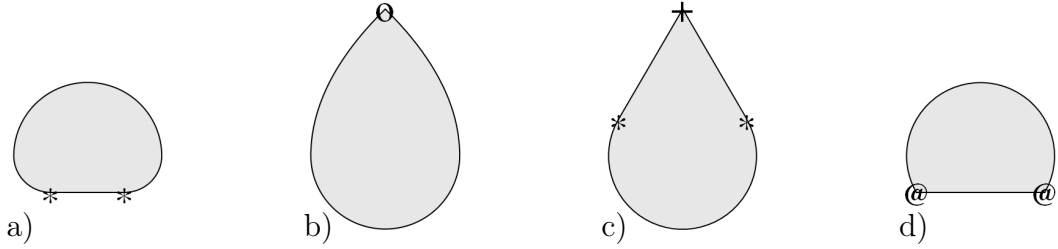
\begin{Exa}
\label{ex:2}
The convex body in Figure~\ref{fig:dual_pics} a) appears at
$a=b=\tfrac 1{2}$ in a family of convex bodies $K\subset\mathbb R^2$ 
defined for $a,b>0$. The dual convex body $K^*$ is depicted in 
Figure~\ref{fig:dual_pics} b) for $a=b=\tfrac 1{2}$. We denote
$u(\alpha):=\s{\cos(\alpha)\\\sin(\alpha)}$ for $\alpha\in\mathbb R$.
The boundary of $K$ consists of the segment between $\s{-a\\-b}$ and 
$\s{a\\-b}$, one half arc and two quarter arcs
\[\textstyle
c:\;[0,2\pi)\;\to\;\mathbb R^2\,,\quad
\alpha\;\mapsto\;\left\{\begin{array}{ll}
(a+b)u(\alpha) & \text{for } 0\leq\alpha<\pi\,,\\
\s{-a\\0}+b u(\alpha) & 
\text{for } \pi\leq\alpha<\tfrac{3}{2}\pi\,,\\
\s{a\\0}+b u(\alpha) & 
\text{for } \tfrac{3}{2}\pi\leq\alpha<2\pi\,.
\end{array}\right.
\]
By Construction~\ref{con:support_radial} we have for $\alpha\in[0,2\pi)$
\[\textstyle
\rho_{K^*}(u(\alpha))
\;=\;
\left\{\begin{array}{ll}
(a\cos(\alpha)+b)^{-1} & \text{for } 0\leq\alpha<\tfrac{\pi}{2}\,,\\
(-a\cos(\alpha)+b)^{-1} & \text{for } \tfrac{\pi}{2}\leq\alpha<\pi\,,\\
(a+b)^{-1} & 
\text{for } \pi\leq\alpha\leq 2\pi\,.
\end{array}\right.
\]
\end{Exa}
\par
We define faces and conjugate faces and we prove technical assertions
for \S\ref{sec:2D}.
\begin{Def}
\label{def:faces}
\begin{enumerate}[1.]
\item
A {\it face} of a convex subset $C\subset\mathbb R^n$ is a convex
subset $F\subset C$ such that $x,y,z\in C$, $y\in F$ and 
$y\in\,]x,z[\,:=\{(1-\lambda)x+\lambda z\mid 0<\lambda<1\}$ implies
$x,z\in F$. 
\item
If $u\in\mathbb R^n$ is non-zero then we define
$H_C(u):=\{x\in\mathbb R^n\mid\langle x,u\rangle=h_C(u)\}$.
If $C\cap H_C(u)\neq\emptyset$ then $H_C(u)$ is an affine hyperplane
called {\it supporting hyperplane} and $C\cap H_C(u)$ is an 
{\it exposed face} of $C$. By definition $\emptyset$ and $C$ are
exposed faces of $C$. A face which is not an exposed face is called a
{\it non-exposed face}.
\item
If $\{x\}$ is a face of $C$ for $x\in C$ then $x$ is an
{\it extremal point}. In the following we will identify extremal points 
with their faces. If the extremal point $x\in C$ is an exposed face then
$x$ is an {\it exposed point}, otherwise $x$ is a {\it non-exposed} point.
\item
The {\it conjugate face} $\mathcal C(F)$ of a subset $F\subset K$ is a
subset of the polar convex body:
\begin{equation}
\label{eq:conjugate_face}\textstyle
\mathcal C(F)
\;=\;
\mathcal C_K(F)
\;:=\;
\{v\in K^\circ\mid\langle v,u\rangle=1\;\forall u\in F\}\,.
\end{equation}
\end{enumerate}
\end{Def}
\begin{Rem}
Exposed faces of a convex subset $C\subset\mathbb R^n$ are faces of $C$,
see e.g.\ \S18 in \cite{Rockafellar}. It is a common practice to use the
conjugate face mapping $\mathcal C$ without reference to the convex body 
$K$ and write e.g.\ $\mathcal C^2(F)$, see \S2.2 in \cite{Schneider}. 
\end{Rem}
\par
We denote $H^\pm:=\{(x,y)^T\in\mathbb R^2\mid\pm y\geq 0\}$.
\begin{Lem}
\label{lem:ind_q}
Let $L\subset\mathbb R^2$ be a convex body.
\begin{enumerate}[1.]
\item
A point $x\in L\setminus H^\mp$ is an extremal point of 
$L$ if and only if $x$ is an extremal point of $L\cap H^\pm$.
\item
Let $\rho_L(\pm\s{1\\0})=h_L(\pm\s{1\\0})$, i.e.\ $L$ has maximal 
$x$-extension on the $x$-axis.
\begin{enumerate}[a)]
\item
The support functions satisfy 
$h_L|_{H^{\pm}}=h_{L\cap H^\pm}|_{H^{\pm}}$.
\item
For every $u\in\mathbb R^2\setminus H^\mp$ the supporting hyperplanes
satisfy $H_L(u)=H_{L\cap H^\pm}(u)$.
\item
If $p\in L\setminus H^\mp$ is an exposed point of $L$ then 
there is $u\in\mathbb R^2\setminus H^\mp$ such that 
$\{p\}=L\cap H_L(u)$. 
\item
If $0\in{\rm int}(L)$ and $F\subset L$ such that 
$F\setminus H^\mp\neq\emptyset$, then $\mathcal C_L(F)\subset H^\pm$.
\end{enumerate}
\item
For $i=1,2$ let $L_i\subset\mathbb R^2$ be a convex body
with $0\in{\rm int}(L_i)$ and let $c_\pm>0$ such that
$\rho_{L_i}(\pm\s{1\\0})=h_{L_i}(\pm\s{1\\0})=c_\pm$.
Then $L:=(L_1\cap H^+)\cup(L_2\cap H^-)$ is a convex body with
$0\in{\rm int}(L)$ and $L^*=(L_2^*\cap H^+)\cup(L_1^*\cap H^-)$.
\end{enumerate}
\end{Lem}
{\em Proof:\/}
The proof of 1 and 2 is written for $(\pm,\mp)=(+,-)$,  
$(\pm,\mp)=(-,+)$ is analogous.
To show part~1 let $x\in L\setminus H^-$. If $x$ is an extremal point
of $L$ then it is trivially an extremal point of $L\cap H^+$.
Conversely let $x$ be an extremal point of $L\cap H^+$ and let
$y,z\in L$ with $x\in\,]y,z[$. The case $y,z\in L\setminus H^+$ is 
impossible since $x\in H^+$. If $y,z\in L\cap H^+$ then $y=z=x$
follows as desired. Finally, if $y\in L\cap H^+$ and 
$z\in L\setminus H^+$, then $]y,z[$ intersects the $x$-axis in a
point $p\neq x$. Then as before $y=p=x$ and this implies $z=x$.
\par
To prove part 2 a) we show for $u=(u_x,u_y)\in H^+$ that
$\langle\cdot,u\rangle$ is maximized on $L$ at a point in $L\cap H^+$.
Let $p=(p_x,p_y)\in L\cap H^-$. Assuming $\pm u_x\geq 0$ we show 
$\langle p,u\rangle\leq\langle\pm\s{1\\0}\rho_L(\pm\s{1\\0}),u\rangle$.
Since $L$ satisfies $h_L(\pm\s{1\\0})=\rho_L(\pm\s{1\\0})$ we have
\[\textstyle
\pm p_x
\;=\;
\langle p,\pm\s{1\\0}\rangle
\;\leq\;
h_L(\pm\s{1\\0})
\;=\;
\rho_L(\pm\s{1\\0})
\]
and
\begin{align}
\label{eq:inequality}
\langle p,u\rangle
&\;=\;p_xu_x+p_yu_y
\;\leq\;
p_xu_x
\;=\;
(\pm p_x)(\pm u_x)
\;\leq\;
\rho_L(\pm\s{1\\0})(\pm u_x)\\\nonumber
&\;=\;
\langle\pm\s{1\\0}\rho_L(\pm\s{1\\0}),u\rangle\,.
\end{align}
\par
The assertion 2 b) holds because $u\not\in H^-$ and $p\not\in H^+$
imply $p_yu_y<0$ and then a strict inequality follows in
(\ref{eq:inequality}). 
\par
We show 2 c). Since $p$ is an exposed point of $L$, there exists a 
non-zero vector $u$ with $\{p\}=L\cap H_L(u)$. By contradiction
we show $u\not\in H^-$. By 2 a) there is
point $q\in L\cap H^-$ that lies on the hyperplane $H_L(\pm\s{1\\0})$. 
Since $p\not\in H^-$ the vector $u$ is not aligned with the $x$-axis.
If we assume $u\not\in H^+$ then 2 b) shows $p\not\in H_L(u)$.
\par
For 2 d) we show that $F\setminus H^-\neq\emptyset$ implies
$\mathcal C(F)\subset H^+$ by proving $p\not\in\mathcal C(F)$ 
for every $p=(p_x,p_y)^T$ in $L^\circ\setminus H^+$. 
We have $p_y<0$ and there exists $u=(u_x,u_y)^T\in F$
such that $u_y>0$. Since $0\in{\rm int}(L)$ the polar $L^\circ$
is a convex body and by (\ref{eq:rad_supp}) it 
satisfies $\rho_{L^\circ}(\pm\s{1\\0})=h_{L^\circ}(\pm\s{1\\0})$.
Assuming $\pm u_x\geq 0$ the strict inequality $\langle p,u\rangle
<\langle\pm\s{1\\0}\rho_{L^\circ}(\pm\s{1\\0}),u\rangle$ follows 
from (\ref{eq:inequality}) with $L$ replaced by $L^\circ$. Since
$\pm\s{1\\0}\rho_{L^\circ}(\pm\s{1\\0})\in L^\circ$ and $u\in L$
we have
$\langle\pm\s{1\\0}\rho_{L^\circ}(\pm\s{1\\0}),u\rangle\leq 1$ 
hence $\langle p,u\rangle<1$ shows $p\not\in\mathcal C(F)$.
\par
We show part 3.
Clearly $L$ is compact and $0\in{\rm int}(L)$. To show convexity let 
$x,y\in L$ and 
$[x,y]:=\{(1-\lambda)x+\lambda y\mid\lambda\in[0,1]\}$. If
$x,y\in H^\pm$ then $[x,y]\subset L$ by convexity of $L_1$ and $L_2$
Otherwise $[x,y]$ intersects the $x$-axis in a point $p$ and the
pairs $\{x,p\}$ and $\{p,y\}$ satisfy the previous assumption. 
Using (\ref{eq:rad_supp}) and 2 b) we have 
\[\textstyle
\rho_{L^*}(u)
\;=\;h_L(-u)^{-1}
\;=\;h_{L_2}(-u)^{-1}
\;=\;\rho_{L_2^*}(u)
\]
for all $u\in H^+$. Similarly for $u\in H^-$ we have 
$\rho_{L^*}(u)=\rho_{L_1^*}(u)$. 
\hspace*{\fill}$\Box$\\
\par
The following construction of planar self-dual convex bodies joins half
of a convex body with half of its dual convex body. By part b) the 
construction is general.
\begin{Con}
\label{con:assemble_two_half}
\begin{enumerate}[a)]
\item
Let $K\subset\mathbb R^2$ satisfy 
$\rho_K(\pm\s{1\\0})=h_K(\pm\s{1\\0})=e^{\pm\lambda}$ for some 
$\lambda\in\mathbb R$. Then $(K\cap H^+)\cup(K^*\cap H^-)$ is a self-dual 
convex body. 
\item
For every planar self-dual convex body $K$ exists a rotation 
$\psi\in SO(2)$ such that $\psi(K)$ satisfies the assumptions in a).
\end{enumerate}
\end{Con}
{\em Proof:\/}
Assertion a) follows from (\ref{eq:rad_supp}) and
Lemma~\ref{lem:ind_q}.3 applied to the convex bodies $L_1:=K$ and 
$L_2:=K^*$. 
To show b) let $u$ be an element of $K$ with maximal norm $|u|$ in $K$
and put $v:=\tfrac u{|u|}$. Then $\rho_K(v)=|u|$ and
$h_K(v)=\max_{w\in K}\langle w,v\rangle\leq\max_{w\in K}|w||v|=|u|$ by 
the Cauchy-Schwarz inequality. On the other hand,
$h_K(v)\geq\langle u,v\rangle=|u|$ shows $\rho_K(v)=h_K(v)$. Since $K$
is self-dual we get from (\ref{eq:rad_supp}) and with $(K^*)^*=K$
\[\textstyle
\rho_K(-v)
\;=\;\rho_{K^*}(-v)
\;=\;h_K(v)^{-1}
\;=\;\rho_K(v)^{-1}
\;=\;h_{K^*}(-v)
\;=\;h_{K}(-v),
\]
that is $\rho_K(\pm v)=h_K(\pm v)=e^{\pm\lambda}$ for some 
$\lambda\in\mathbb R$. For all $\psi\in SO(2)$ and $v\in\mathbb R^2$ the 
equalities $\rho_{\psi(K)}(\psi(v))=\rho_K(v)$ and 
$h_{\psi(K)}(\psi(v))=h_K(v)$ hold. The choice of $\psi$ such that
$\psi(v)=\s{1\\0}$ completes the proof. 
\hspace*{\fill}$\Box$\\
\begin{figure}
\centerline{%
\begin{picture}(2.14,3.2)
\put(0,0){\includegraphics[height=3.2cm, bb=0 0 267 400, 
clip=]{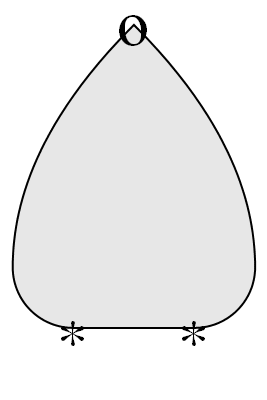}}
\put(0,0){a)}
\end{picture}
\hspace{1.0cm}
\begin{picture}(2.14,3.2)
\put(0,0){\includegraphics[height=3.2cm, bb=0 0 267 400, 
clip=]{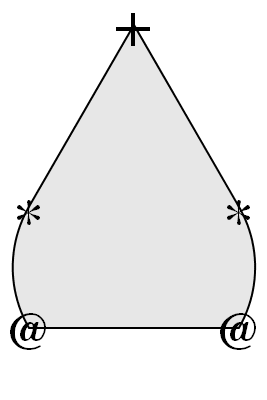}}
\put(0,0){b)}
\end{picture}
\hspace{1.0cm}
\begin{picture}(2.33,3.2)
\put(0,0){\includegraphics[height=3.2cm, bb=0 0 291 400, 
clip=]{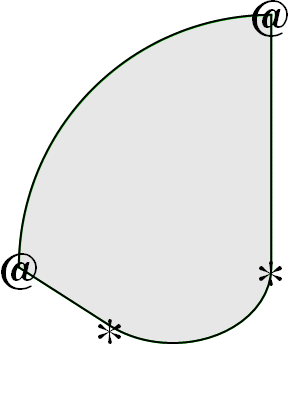}}
\put(0,0){c)}
\end{picture}
\hspace{1.5cm}
\begin{picture}(3.2,3.2)
\put(0,0){\includegraphics[height=3.2cm, bb=0 0 401 400, 
clip=]{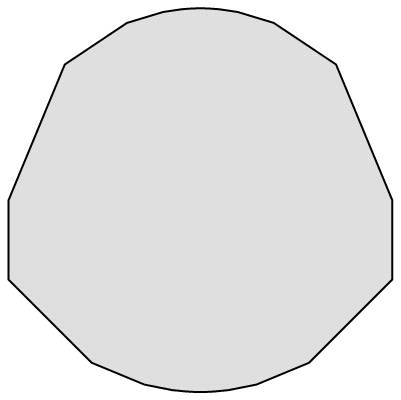}}
\put(0,0){d)}
\end{picture}}
\caption{\label{fig:self_pics} The depicted convex sets are self-dual. 
Markings are explained in Figure~\ref{fig:dual_pics}. Drawings a)--c)
have equal scaling and their origin is aligned vertically.}
\end{figure}
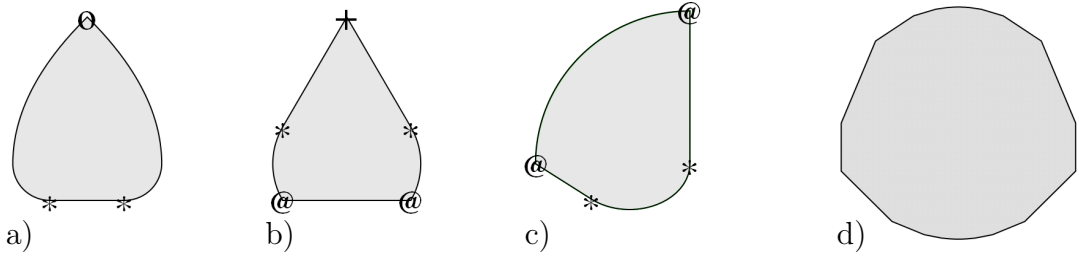
\par
The self-dual convex body a) resp.\ b) in Figure~\ref{fig:self_pics} is 
generated by Construction~\ref{con:assemble_two_half} from the convex 
body b) resp.\ c) and its dual convex body a) resp.\ d) in 
Figure~\ref{fig:dual_pics}. We consider a less symmetric example.
\begin{Exa}
Let $a>0$ and $K\subset\mathbb R^2$ have the upper part
$K\cap H^+$ defined as the convex hull of $\s{a\\0}$ and of the quarter 
arc consisting of all points $\s{a\\0}+\tfrac{a^2+1}{a}u(\alpha)$ for 
$\alpha\in[\tfrac{\pi}{2},\pi]$. Without specifying the lower 
part of $K$, Construction~\ref{con:support_radial}, 
Lemma~\ref{lem:ind_q}.2 a) and Construction~\ref{con:assemble_two_half} 
provide a self-dual convex body $X:=(K\cap H^+)\cup(K^*\cap H^-)$ with
radial function 
\begin{align*}
\rho_X(u(\alpha))
&\;=\;
\rho_{K^*}(u(\alpha))
\;=\;
h_K(-u(\alpha))^{-1}\\
&\;=\;
\left\{\begin{array}{ll}
-a(a^2(\cos(\alpha) + \sin(\alpha)) + \sin(\alpha))^{-1}
 & \text{for } \pi\leq\alpha<\tfrac 3{2}\pi,\\
a (a^2 (1 - \cos(\alpha)) + 1)^{-1} & 
\text{for } \tfrac 3{2}\pi\leq\alpha\leq 2\pi.
\end{array}\right.
\end{align*}
For $a=\tfrac 4{3}$ the self-dual convex body $X$ is shown in 
Figure~\ref{fig:self_pics} c). 
\end{Exa}
%
%
%
%
%
%
\section{A Galois connection}
\label{sec:galois}
\par
We define a Galois connection between touching cones and faces of an 
arbitrary convex subset $C\subset\mathbb R^n$ which has not cardinality 
one. We will study two lattices of faces and two lattices of cones 
associated to $C$. We refer to \cite{Birkhoff} for general lattice
theory and to ~\cite{Barker_Lattice,Loewy,Weis_cones} for the lattice 
theory of convex sets.
\par
The {\it normal cone} at a point $x\in C$ is the set $N(x)$ of all vectors 
$u\in\mathbb R^n$ such that $\langle u,y-x\rangle\leq 0$ holds for all
$y\in C$, i.e.\ $u$ does not make an acute angle with $y-x$ for any
$y\in C$. The whole space $\mathbb R^n$ is a normal cone by definition.
A {\it touching cone} is any non-empty face of any normal cone of $C$. 
(Touching cones were first introduced in \cite{Schneider} by a different
but equivalent definition).
\par
The set of faces, exposed faces, touching cones resp.\ normal cones of
$C$ is denoted by $\mathcal F_C$, $\mathcal E_C$, $\mathcal T_C$
resp.\ $\mathcal N_C$. We have 
\[\textstyle
\mathcal E_C\;\subset\;\mathcal F_C
\quad\text{and}\quad
\mathcal N_C\;\subset\;\mathcal T_C\,.
\]
Each of these sets is a poset ordered by inclusion and a complete lattice
of finite length where the infimum is the intersection, see e.g.\
\S1.1 and \S1.2 in \cite{Weis_cones}. We write these lattices in the form
\begin{align*}
(\mathcal F_C,\leq_{\mathcal F},\vee_{\mathcal F},
\wedge_{\mathcal F})\,,\quad
&(\mathcal E_C,\leq_{\mathcal E},\vee_{\mathcal E},
\wedge_{\mathcal E})\,,\\
(\mathcal T_C,\leq_{\mathcal T},\vee_{\mathcal T},
\wedge_{\mathcal T})\,,\quad
&(\mathcal N_C,\leq_{\mathcal N},\vee_{\mathcal N},
\wedge_{\mathcal N})\,.
\end{align*}
The infimum (supremum) of a subset $S\subset\mathcal F_C$ is denoted by
$\bigwedge_{\mathcal F}S$ ($\bigvee_{\mathcal F}S$), the analogue
notation is used for other lattices.
\par
We consider a mapping $\theta:L\to M$ between two lattices
$(L, \leq_L, \vee_L, \wedge_L)$ and $(M, \leq_M, \vee_M, \wedge_M)$.
The mapping $\theta$ is 
\begin{align*}
\text{{\it isotone} if }
& \quad x\leq_L y 
\;\implies\;
\theta(x)\leq_M\theta(y)\,,
& (x,y\in L)\\
\text{{\it antitone} if }
& \quad x\leq_L y 
\;\implies\;
\theta(x)\geq_M\theta(y)\,,\\
\text{a {\it join-morphism} if }
& \quad \theta(x\vee_L y)
\;=\;
\theta(x)\vee_M\theta(y)\,,\\
\text{a {\it meet-morphism} if }
& \quad \theta(x\wedge_L y)
\;=\;
\theta(x)\wedge_M\theta(y)\,,\\
\text{a {\it dual join-morphism} if }
& \quad \theta (x\vee_L y)
\;=\;
\theta(x)\wedge_M\theta(y)\\
\text{and a {\it dual meet-morphism} if }
& \quad \theta(x\wedge_L y)
\;=\;
\theta(x)\vee_M\theta(y)\,.
\end{align*}
Join- and meet-morphisms are isotone, see \cite{Birkhoff} Chap.\ 
II.3. Hence dual join- and meet-morphisms are antitone. A
{\it (dual) lattice-morphism} is a (dual) meet-morphism which is 
also a (dual) join-morphism.
\par
The {\it relative interior} of $C$, denoted by ${\rm ri}(C)$ is the
interior of $C$ in the topology of the affine hull ${\rm aff}(C)$ of $C$.
If $C\neq\emptyset$ then the translation vector space of the affine hull
of $C$ is denoted by ${\rm lin}(C):={\rm aff}(C)-{\rm aff}(C)$.
\begin{Def}
\label{def:varphiN}
\begin{enumerate}[1.]
\item
To every touching cone we associate an exposed face
\[\textstyle
\Phi:\;\mathcal T_C\;\longrightarrow\;\mathcal E_C\,.
\]
We put $\Phi({\rm lin}(C)^\perp):=C$,
$\Phi(\mathbb R^n):=\emptyset$ and for touching cones
$T\in\mathcal T_C\setminus\{{\rm lin}(C)^\perp,\mathbb R^n\}$
Lemma 7.2 a) in \cite{Weis_cones} shows that the exposed face
\[\textstyle
\Phi(T)\;:=\;C\cap H_C(u)
\]
is well-defined for an arbitrary non-zero vector 
$u$ in the relative interior of $T$.
\item
To every face we associate a normal cone
\[\textstyle
\Psi:\;\mathcal F_C\;\longrightarrow\;\mathcal N_C\,.
\]
We put $\Psi(\emptyset):=\mathbb R^n$. For faces
$F\in\mathcal F_C\setminus\{\emptyset\}$ a normal cone is well-defined 
by the arguments provided in Definition 4.3 in \cite{Weis_cones}: We put
\[\textstyle
\Psi(F)\;:=\;N(x)
\]
for an arbitrary point $x$ in the relative interior of $F$.
\end{enumerate}
\end{Def}
\begin{Rem}
\label{rem:varphiNanti}
\begin{enumerate}[1.]
\item
The map $\Phi:\mathcal T_C\to\mathcal E_C$ is antitone, this
follows from an intersection representation: If $T\in\mathcal T_C$
is a touching cone and $T\neq{\rm lin}(C)^\perp,\mathbb R^n$ then by
Lemma 7.2 a) in \cite{Weis_cones} we have
\[\textstyle
\Phi(T)
\;=\;
\bigcap_{u\in T\setminus\{0\}}(C\cap H_C(u))\,.
\]
\item
That $\Psi:\mathcal F_C\to\mathcal N_C$ is antitone is discussed in
the paragraph following Definition 4.3 in \cite{Weis_cones}. For the
sake of completeness we notice for faces $F\in\mathcal F_C$,
$F\neq\emptyset$:
\[\textstyle
\Psi(F)
\;=\;
\bigcap_{x\in F}N(x)\,.
\]
This follows from the inclusion $N(y)\subset N(x)$ valid for all
$y$ in the relative interior ${\rm ri}F$ and $x\in F$, see (15)(ii)
in \cite{Weis_cones}, while $\Psi(F)=N(y)$ holds by
Definition~\ref{def:varphiN}.2.
\item
It is proved in Proposition 4.7 in \cite{Weis_cones} that the
restrictions $\Phi|_{\mathcal N_C}$ and $\Psi|_{\mathcal E_C}$ are
dual lattice isomorphisms, inverse to each other. The diagram 
\begin{equation}
\label{eq:l_anti_normalcone}\textstyle
\xymatrix{%
\mathcal N_C\quad \ar@{^{(}->>}[rr]<0.35ex>^{\Phi} 
 & & \quad\mathcal E_C \ar@{^{(}->>}[ll]<0.35ex>^{\Psi}
}
\end{equation}
commutes. 
\end{enumerate}
\end{Rem}
\par
To study $\Phi$ and $\Psi$ we use the concepts of closure operation and
of Galois connection, see e.g.\ \S{}V.1 and \S{}V.8 in \cite{Birkhoff}.
\begin{Def}
\label{def:closureop}
\begin{enumerate}[1.]
\item
A {\it closure operation} on a set $I$ is an operator $X\to{\rm cl}(X)$
on the subsets of $I$ such that for all $X,Y\subset I$ we have\\[3mm]
\centerline{%
\begin{tabular}{lcr}
$X\;\subset\;{\rm cl}(X)$ & \qquad & (Extensive)\\
${\rm cl}(X)\;=\;{\rm cl}({\rm cl}(X))$ & & (Idempotent)\\
If $X\subset Y$, then ${\rm cl}(X)\subset{\rm cl}(Y)$ & &
(Isotone)
\end{tabular}}\\[3mm]
Subsets $X\subset I$ with $X={\rm cl}(X)$ are called {\it closed sets}
with respect to ${\rm cl}$.
\item
Let $(L,\leq_L)$ and $(M,\leq_M)$ be any posets and let $\theta:L\to M$, 
$\phi:M\to L$ be maps such that for all $l_1,l_2\in L$ and 
$m_1,m_2\in M$ we have\\[3mm]
\centerline{%
\begin{tabular}{c}
$l_1\;\leq_L\;l_2$
\quad implies \quad
$\theta(l_1)\;\geq_M\;\theta(l_2)$\,,\\
$m_1\;\leq_M\;m_2$
\quad implies \quad
$\phi(m_1)\;\geq_L\;\phi(m_2)$\,,\\
$l_1\;\leq_L\;\phi(\theta(l_1))$
\quad and \quad
$m_1\;\leq_M\;\theta(\phi(m_1))$\,.
\end{tabular}}\\[3mm]
Then $\theta$ and $\phi$ are said to define a {\it Galois connection}
between $L$ and $M$.
\item
We define the {\it normal closure} as the operation on touching cones
\begin{equation}
\label{def:clN}\textstyle
{\rm cl}_{\mathcal N}:\;\mathcal T_C\;\longrightarrow\;
\mathcal N_C\,,\quad
T\;\longmapsto\;\bigcap_{\substack{N\in\mathcal N_C\\T\subset N}}N
\end{equation}
and the {\it exposed closure} as the operation on faces
\begin{equation}
\label{def:clvarphi}\textstyle
{\rm cl}_{\mathcal E}:\;\mathcal F_C\;\longrightarrow\;
\mathcal E_C\,,\quad
F\;\longmapsto\;\bigcap_{\substack{G\in\mathcal E_C\\F\subset G}}G\,.
\end{equation}
\end{enumerate}
\end{Def}
\par
Since $\mathcal N_C$ and $\mathcal E_C$ are complete lattices with 
the intersection as the infimum, the normal closure and the exposed 
closure are closure operations in the sense of 
Definition~\ref{def:closureop}.1. The closed sets of
${\rm cl}_{\mathcal N}$ are the normal cones and the closed sets of
${\rm cl}_{\mathcal E}$ are the exposed faces. 
\par
These closures can equivalently be defined by the mappings $\Phi$ and
$\Psi$ between touching cones and faces.
\begin{Lem}
\label{lem:Nvarphiclosure}
\begin{enumerate}[1.]
\item
Every touching cone $T\in\mathcal T_C$ has normal closure
${\rm cl}_{\mathcal N}(T)=\Psi\circ\Phi(T)$. In particular
$T\leq_{\mathcal T}\Psi\circ\Phi(T)$ holds.
\item
Every face $F\in\mathcal F_C$ has exposed closure
${\rm cl}_{\mathcal E}(F)=\Phi\circ \Psi(F)$. In particular
$F\leq_{\mathcal F}\Phi\circ \Psi(F)$ holds.
\end{enumerate}
\end{Lem}
{\em Proof:\/}
We prove part 1.
Let $T\in\mathcal T_C$ be a touching cone and let $N\in\mathcal N_C$
be a normal cone. By Remark~\ref{rem:varphiNanti}.1 and~2 the maps 
$\Phi:\mathcal T_C\to\mathcal E_C$ and $\Psi:\mathcal F_C\to\mathcal N_C$
are antitone hence the composition $\Psi\circ\Phi$ is isotone. Its 
restriction $\Psi\circ\Phi|_{\mathcal N_C}$ is the identity mapping by 
(\ref{eq:l_anti_normalcone}) hence
\[\textstyle
T\subset N
\;\implies\;
\Psi\circ\Phi(T)\subset\Psi\circ\Phi(N)=N\,.
\]
This implication has two consequences. Firstly, the inclusion
$\Psi\circ\Phi(T)\subset{\rm cl}_{\mathcal N}(T)$ into the normal closure 
(\ref{def:clN}) follows. Secondly we have
\[\textstyle
{\rm cl}_{\mathcal N}(T)
\;\stackrel{\text{def.}}{=}\;
\bigcap_{\substack{N\in\mathcal N_C\\T\subset N}}N
\;\subset\;
\bigcap_{\substack{N\in\mathcal N_C\\\Psi\circ\Phi(T)\subset N}}N
\;=\;
\Psi\circ\Phi(T)
\]
where the last equality holds because $\Psi\circ\Phi(T)\in\mathcal N_C$.
This shows $\Psi\circ\Phi(T)={\rm cl}_{\mathcal N}(T)$. The inclusion 
$T\subset{\rm cl}_{\mathcal N}(T)$ is obvious. The proof of part 2 is
analogous.
\hspace*{\fill}$\Box$\\
\par
Their link to the closure operations enables us to analyze
$\Phi$ and $\Psi$. 
\begin{Lem}
\label{lem:dual_join}
\begin{enumerate}[1.]
\item
The assignment of exposed faces to touching cones
$\Phi:\mathcal T_C\to\mathcal E_C$ is a dual join-morphism.
For all touching cones $T,U\in\mathcal T_C$ we have
$\Phi(T\wedge_{\mathcal T}U)\geq_{\mathcal E}
\Phi(T)\vee_{\mathcal E}\Phi(U)$ and
$\Phi(T)=\Phi({\rm cl}_{\mathcal N}(T))$.
\item
The assignment of normal cones to faces
$\Psi:\mathcal F_C\to\mathcal N_C$ is a dual join-morphism. For all
faces $F,G\in\mathcal F_C$ we have 
$\Psi(F\wedge_{\mathcal F}G)\geq_{\mathcal N}
\Psi(F)\vee_{\mathcal N}\Psi(G)$ and 
$\Psi(F)=\Psi({\rm cl}_{\mathcal E}(F))$.
\end{enumerate}
\end{Lem}
{\em Proof:\/}
We prove part 1 in five steps. As $\Phi$ is antitone by
Remark~\ref{rem:varphiNanti}, for all touching cones $T,U\in\mathcal T_C$
follows
\begin{align*}
T\wedge_{\mathcal T} U\;\leq_{\mathcal T}\;T,U
& \quad\implies\quad
\Phi(T\wedge_{\mathcal T} U)\;\geq_{\mathcal E}\;\Phi(T),\Phi(U)\\
& \quad\implies\quad
\Phi(T\wedge_{\mathcal T} U)
\;\geq_{\mathcal E}\;\Phi(T)\vee_{\mathcal E}\Phi(U)\quad
\text{and secondly}\\
T\vee_{\mathcal T} U\;\geq_{\mathcal T}\;T,U
& \quad\implies\quad
\Phi(T\vee_{\mathcal T} U)\;\leq_{\mathcal E}\;\Phi(T),\Phi(U)\\
& \quad\implies\quad
\Phi(T\vee_{\mathcal T} U)
\;\leq_{\mathcal E}\;\Phi(T)\wedge_{\mathcal E}\Phi(U)\,.
\end{align*}
Thirdly, by (\ref{eq:l_anti_normalcone}) and
Lemma~\ref{lem:Nvarphiclosure}.1 we have
\[\textstyle
\Phi(T)
\;=\;
\Phi\circ\Psi\circ\Phi(T)
\;=\;
\Phi({\rm cl}_{\mathcal N}(T))\,.
\]
Fourthly, as $\mathcal N_C$ is a complete lattice with the restricted 
partial order from $\mathcal T_C$,
\[\textstyle
T\vee_{\mathcal T}U
\;\leq_{\mathcal T}\;
{\rm cl}_{\mathcal N}(T)\vee_{\mathcal T}{\rm cl}_{\mathcal N}(U)
\;\leq_{\mathcal T}\;
{\rm cl}_{\mathcal N}(T)\vee_{\mathcal N}{\rm cl}_{\mathcal N}(U)\,.
\]
Finally, by step two, step three, since
$\Phi|_{\mathcal N_C}:\mathcal N_C\to\mathcal E_C$ is a dual lattice
isomorphism and by step four we have
\begin{align*}
\Phi(T\vee_{\mathcal T}U)
& \;\leq_{\mathcal E}\;
\Phi(T)\wedge_{\mathcal E}\Phi(U)
\;=\;
\Phi({\rm cl}_{\mathcal N}(T))
\wedge_{\mathcal E}\Phi({\rm cl}_{\mathcal N}(U))\\
& \;=\;
\Phi({\rm cl}_{\mathcal N}(T)\vee_{\mathcal N}{\rm cl}_{\mathcal N}(U))
\;\leq_{\mathcal E}\;
\Phi(T\vee_{\mathcal T}U)\,.
\end{align*}
This completes the proof of part 1, part 2 is analogous.
\hspace*{\fill}$\Box$\\
\begin{Exa} 
\label{exa:counterexample}
The dual join morphisms in Lemma~\ref{lem:dual_join} are no dual lattice
morphisms in general. A counterexample for $\Psi$ is given by the two 
non-exposed faces of the convex set in Figure~\ref{fig:dual_pics} a) or c)
which is also a counterexample against a lattice morphism of the
exposed closure in Lemma~\ref{lem:closures}. The convex bodies 
in Figure~\ref{fig:dual_pics} b) or d) are counterexamples
for $\Phi$ and for the normal closure, because b) is dual to a) and 
d) is dual to c), see Proposition~\ref{pro:diagram}.
\end{Exa}
\par
The closure operations inherit properties from $\Phi$ and $\Psi$.
\begin{Lem}
\label{lem:closures}
\begin{enumerate}[1.]
\item
The normal closure ${\rm cl}_{\mathcal N}:\mathcal T_C\to\mathcal N_C$
is a join morphism such that for all touching cones
$T,U\in\mathcal T_C$ we have
${\rm cl}_{\mathcal N}(T\wedge_{\mathcal T}U)\leq_{\mathcal N}
{\rm cl}_{\mathcal N}(T)\wedge_{\mathcal N}{\rm cl}_{\mathcal N}(T)$.
\item
The exposed closure ${\rm cl}_{\mathcal E}:\mathcal F_C\to\mathcal E_C$
is a join morphism such that for all faces $F,G\in\mathcal F_C$ we have
${\rm cl}_{\mathcal E}(F\wedge_{\mathcal F}G)\leq_{\mathcal E}
{\rm cl}_{\mathcal E}(F)\wedge_{\mathcal E}{\rm cl}_{\mathcal E}(G)$.
\end{enumerate}
\end{Lem}
{\em Proof:\/}
We prove part 1 and choose touching cones $T,U\in\mathcal T_C$. By
Lemma~\ref{lem:Nvarphiclosure}.1, Lemma~\ref{lem:dual_join}.1,
(\ref{eq:l_anti_normalcone}) and Lemma~\ref{lem:Nvarphiclosure}.1 we
have
\begin{align*}
{\rm cl}_{\mathcal N}(T\vee_{\mathcal T}U)
&\;=\;
\Psi\circ\Phi(T\vee_{\mathcal T}U)
\;=\;
\Psi(\Phi(T)\wedge_{\mathcal E}\Phi(U))\\
&\;=\;
\Psi\circ\Phi(T)\vee_{\mathcal N}\Psi\circ\Phi(U)
\;=\;
{\rm cl}_{\mathcal N}(T)\vee_{\mathcal N}{\rm cl}_{\mathcal N}(U)\,.
\end{align*}
The same arguments as above prove
\begin{align*}
{\rm cl}_{\mathcal N}(T\wedge_{\mathcal T}U)
&\;=\;
\Psi\circ\Phi(T\wedge_{\mathcal T}U)
\;\leq_{\mathcal N}\;
\Psi(\Phi(T)\vee_{\mathcal E}\Phi(U))\\
&\;=\;
\Psi\circ\Phi(T)\vee_{\mathcal N}\Psi\circ\Phi(U)
\;=\;
{\rm cl}_{\mathcal N}(T)\vee_{\mathcal N}{\rm cl}_{\mathcal N}(U)\,,
\end{align*}
except the inequality $\leq_{\mathcal N}$ follows because $\Psi$ is 
antitone by Remark~\ref{rem:varphiNanti}.2 and because
Lemma~\ref{lem:dual_join}.1 shows
$\Phi(T\wedge_{\mathcal T}U)\geq_{\mathcal E}
\Phi(T)\vee_{\mathcal E}\Phi(U)$. The proof of part 2 is analogous.
\hspace*{\fill}$\Box$\\
\par
We summarize a part of our results as follows.
\begin{Thm}
\label{thm:galois}
Let $C\subset\mathbb R^n$ be any convex subset of cardinality not one. 
Then the map $\Phi:\mathcal T_C\to\mathcal E_C$ from touching cones to
exposed faces and the map $\Psi:\mathcal F_C\to\mathcal N_C$  from faces 
to normal cones define a Galois connection between the touching cone
lattice $\mathcal T_C$ and the face lattice $\mathcal F_C$.
\end{Thm}
{\em Proof:\/}
This follows from Lemma~\ref{lem:Nvarphiclosure} and 
Lemma~\ref{lem:dual_join}.
\hspace*{\fill}$\Box$\\
\par
We recover the dual lattice isomorphism (\ref{eq:l_anti_normalcone}) 
from an abstract theorem:
\begin{Rem}
\label{ref:abstract_closure}
If $\theta:L\to M$, $\phi:M\to L$ is a Galois connection between complete 
lattices $L$ and $M$, then the maps $\phi\circ\theta$ and $\theta\circ\phi$
are closure operations. Moreover, $\theta$ and $\phi$ restricts to a dual 
lattice isomorphism between the complete lattices of closed sets of 
$\phi\circ\theta$ and $\theta\circ\phi$. This is proved in 
\S{}V.8 in \cite{Birkhoff}.
\end{Rem}
%
%
%
%
\section{Conjugate faces of a convex body}
\label{sec:pairs_of_bodies}
\par
We study conjugate faces of a polar pair of convex bodies in a lattice
theoretic perspective. This pair will be given by the convex body
$K\subset\mathbb R^n$ with $0\in{\rm int}(K)$ and by its polar convex
body $K^\circ\subset\mathbb R^n$ with  $0\in{\rm int}(K^\circ)$. 
\par
In the following we consider the conjugate face map
(\ref{eq:conjugate_face}) in the restriction to the face lattice
$\mathcal F_K$ of $K$,
\[\textstyle
\mathcal C_K:\;\mathcal F_K\;\to\;\mathcal E_{K^\circ}\,.
\]
It is obvious by definition that $\mathcal C_K(\mathcal F_K)$ is included
in the exposed face lattice $\mathcal E_{K^\circ}$. Similarly we consider
the conjugate face map 
$\mathcal C_{K^\circ}:\mathcal F_{K^\circ}\to\mathcal E_K$.
\par
It is well-known that the two conjugate face maps $\mathcal C_K$ and
$\mathcal C_{K^\circ}$ define a Galois connection, see
Definition~\ref{def:closureop}, between the face lattices $\mathcal F_K$
and $\mathcal F_{K^\circ}$. The corresponding closure operations, see 
Remark~\ref{ref:abstract_closure}, are the exposed closure
operations~(\ref{def:clvarphi})
\begin{equation}
\label{eq:galois_conj}\textstyle
\mathcal C_{K^\circ}\circ\mathcal C_K
\;=\;
{\rm cl}_{\mathcal E}|_{\mathcal F_K}
\qquad\text{and}\qquad
\mathcal C_K\circ\mathcal C_{K^\circ}
\;=\;
{\rm cl}_{\mathcal E}|_{\mathcal F_{K^\circ}}\,.
\end{equation}
A proof of these statements is given in
Theorem 2.1.4 in \cite{Schneider}. Equation (\ref{eq:galois_conj}) 
brings the Galois connection (\ref{eq:galois_conj}) in contact with 
the Galois connection in Theorem~\ref{thm:galois}, once for $K$ and once
for $K^\circ$. The latter consists of $\Phi$ assigning exposed faces to
touching cones and $\Psi$ assigning normal cones to faces. The arguments
in this paragraph already integrate all solid and dashed arrows into the
diagram in Proposition~\ref{pro:diagram}.
\par
The dotted arrows in the diagram arise from the positive hull operator.
The {\it positive hull} of $X\subset\mathbb R^n$ is
${\rm pos}(X):=\{\lambda x|\lambda\geq 0, x\in X\}$ unless $X=\emptyset$ 
where ${\rm pos}(\emptyset):=\{0\}$. Lemma 2.2.3 in \cite{Schneider}
proves for faces $F\in\mathcal F_K$ (indeed for non-empty convex subsets
of $K$)
\begin{equation}
\label{eq:conj_norm}\textstyle
\Psi(F)\;=\;{\rm pos}\circ\mathcal C_K(F)\,,
\end{equation}
i.e.\ the normal cone $\Psi(F)$ is the positive hull of the conjugate face.
It follows from (\ref{eq:conj_norm}), (\ref{eq:galois_conj}) and the dual 
lattice isomorphism (\ref{eq:l_anti_normalcone}) that we have a lattice
isomorphism
\begin{equation}
\label{eq:half_iso_pos}
{\rm pos}|_{\mathcal E_{K^\circ}}:\mathcal E_{K^\circ}\to\mathcal N_K\,.
\end{equation}
Theorem 8.3 in \cite{Weis_cones} uses (\ref{eq:half_iso_pos}) and an
elementary analysis of sections of normal cones to prove the lattice
isomorphism
\begin{equation}
\label{eq:l_iso_pos}
{\rm pos}|_{\mathcal F_{K^\circ}}:\mathcal F_{K^\circ}\to\mathcal T_K\,.
\end{equation}
The inverse isomorphism is defined for $T\in\mathcal T_K$ with
$T\neq\mathbb R^n$ by
\begin{equation}
\label{eq:inv_pos}
T\mapsto\partial K^\circ\cap T
\end{equation}
and by $\mathbb R^n\mapsto K^\circ$. Here $\partial K^\circ$ 
denotes the boundary of $K^\circ$.
\begin{Pro}
\label{pro:diagram}
The following diagram commutes. The closure operations 
${\rm cl}_{\mathcal N}$ and ${\rm cl}_{\mathcal E}$ are isotone join
morphisms satisfying 
\[\textstyle
f(a\wedge b)\;\leq\;f(a)\wedge f(b)\qquad\forall a,b\,,
\]
their restrictions to normal cones resp.\ exposed faces is the identity
map.
The mappings $\Phi$, $\Psi$ and the conjugate face maps $\mathcal C_K$
and $\mathcal C_{K^\circ}$ are antitone dual join morphism satisfying
\[\textstyle
f(a\wedge b)\;\geq\;f(a)\vee f(b)\qquad\forall a,b\,,
\]
their restrictions to normal cones resp.\ exposed faces define dual
lattice isomorphisms. The positive hull operator ${\rm pos}$ defines
lattice isomorphisms in the diagram.
\begin{align*}
\xymatrix@=1.3cm@M=3mm{%
%
%
\mathcal T_K \ar@{->>}[rddd]<-1ex>_{{\rm cl}_{\mathcal N}} 
\ar@{->>}[rd]^{\Phi}
 &
 &
 &
 & \mathcal F_{K^\circ} \ar@{_{(}.>>}[llll]<-1ex>_{\rm pos} 
 \ar@{->>}[ld]_{\Psi}
 \ar@{->>}[lddd]<+1ex>^{{\rm cl}_{\mathcal E}}
 \ar@{-->>}[llld]<-1ex>_{\mathcal C_{K^\circ}}
 \\
%
%
 & \mathcal E_K \ar@{^{(}.>>}[rr]^{\rm pos} 
 \ar@{_{(}->>}[dd]_{\Psi}
 \ar@{^{(}-->>}[rrdd]<+0.4ex>^{\mathcal C_K}
 &
 & \mathcal N_{K^\circ} \ar@{_{(}->>}[dd]<-1ex>_{\Phi}\\
%
%
\\
%
%
 & \mathcal N_K \ar@{_{(}->>}[uu]<-1ex>_{\Phi}
 &
 & \mathcal E_{K^\circ} \ar@{^{(}.>>}[ll]^{\rm pos}
 \ar@{_{(}->>}[uu]_{\Psi}
 \ar@{^{(}-->>}[uull]<+0.4ex>^{\mathcal C_{K^\circ}}\\
%
%
 \mathcal F_K \ar@{->>}[ruuu]<+1ex>^{{\rm cl}_{\mathcal E}}
 \ar@{_{(}.>>}[rrrr]<-1ex>_{\rm pos} \ar@{->>}[ru]_{\Psi}
 \ar@{-->>}[rrru]<-1ex>_{\mathcal C_K}
 &
 &
 &
 & \mathcal T_{K^\circ} \ar@{->>}[luuu]<-1ex>_{{\rm cl}_{\mathcal N}}
 \ar@{->>}[lu]^{\Phi}
}
\end{align*}
\end{Pro}
{\em Proof:\/}
The commuting diagram was introduced in the above discussion except we 
have to show the equality of functions 
${\rm cl}_{\mathcal N}\circ{\rm pos}={\rm pos}\circ{\rm cl}_{\mathcal E}$
on the domain of the two face lattices $\mathcal F_K$ or
$\mathcal F_{K^\circ}$. We will carry out the proof for
$F\in\mathcal F_{K^\circ}$, the proof for faces of $K$ is analogous. 
By (\ref{def:clN}), (\ref{eq:l_iso_pos}), (\ref{eq:half_iso_pos}), 
(\ref{eq:half_iso_pos}) and (\ref{def:clvarphi}) we have
\begin{align*}
{\rm cl}_{\mathcal N}\circ{\rm pos}(F)
& \;=\; \textstyle\bigwedge_{\mathcal N}
\{N\in\mathcal N_K\mid{\rm pos}(F)\leq_{\mathcal T}N\}\\
& \;=\; \textstyle{\rm pos}\circ{\rm pos}^{-1}\left(\bigwedge_{\mathcal N}
\{N\in\mathcal N_K\mid F\leq_{\mathcal F}{\rm pos}^{-1}(N)\}\right)\\
& \;=\; \textstyle{\rm pos}\left(\bigwedge_{\mathcal E}
\{{\rm pos}^{-1}(N)\in\mathcal E_{K^\circ}\mid
F\leq_{\mathcal F}{\rm pos}^{-1}(N)\}\right)\\
& \;=\; \textstyle{\rm pos}\left(\bigwedge_{\mathcal E}
\{G\in\mathcal E_{K^\circ}\mid
F\leq_{\mathcal F}G\}\right)\\
& \;=\; \textstyle{\rm pos}\circ{\rm cl}_{\mathcal E}(F)\,.
\end{align*}
\par
Lemma~\ref{lem:closures} shows that the closure operations have the
claimed properties. This is shown for $\Phi$ and $\Psi$ in 
Lemma~\ref{lem:dual_join}. Since the conjugate face map
$\mathcal C_K=\Phi\circ{\rm pos}$ is a composition of $\Phi$ with
the positive hull lattice isomorphism (\ref{eq:l_iso_pos}) it has 
the claimed properties. The argument for $\mathcal C_{K^\circ}$ is 
analogous.
\hspace*{\fill}$\Box$\\
\begin{Rem}
\label{rem:thm}
\begin{enumerate}[1.]
\item
The convex bodies in Figure~\ref{fig:dual_pics} a) or c) show
that the conjugate face map is not a dual lattice morphism
(see Example~\ref{exa:counterexample} for the other mappings.)
Equality conditions of a dual join morphism in the inequality 
\[\textstyle
f(a\wedge b)\;\geq\;f(a)\vee f(b)\qquad\forall a,b
\]
were studied in \cite{Barker_Dual} for face lattices of closed convex
cones in relation to modularity of face lattices.
\item
Although the exposed face lattices $\mathcal E_K$ and
$\mathcal E_{K^\circ}$ are dually isomorphic by the conjugate face map,
the face lattices $\mathcal F_K$ and $\mathcal F_{K^\circ}$ are not
dually isomorphic in general. Examples are the dual pairs of convex
bodies in Figure~\ref{fig:dual_pics}.
\end{enumerate}
\end{Rem}
\par
We notice two restricted isomorphisms of the conjugate face map.
For their discussion we introduce further concepts. We call a non-empty face 
$F\in\mathcal F_K$ {\it singular} if its normal cone has dimension at
least two, ${\rm dim}\,\Psi(F)\geq 2$. A non-empty face $F$ is a
{\it corner} of $K$ if ${\rm dim}\,\Psi(F)=n$. A face $F$ of $K$ is
a {\it facet} if ${\rm codim}(F)=1$. Finally, we call a point $x\in K$ 
{\it smooth} if its normal cone has dimension one, ${\rm dim}\,N(x)=1$.
\begin{Cor}
\label{cor:smooth_exp}
The conjugate face 
$\mathcal C_K:\mathcal F_K\to\mathcal E_{K^\circ}$ restricts to a 
bijection 
\[\textstyle
\{\text{ smooth exposed points of }K\}
\;\to\;
\{\text{  smooth exposed points of }K^\circ\}\,.
\]
\end{Cor}
{\em Proof:\/}
The bijection is immediate from Proposition~\ref{pro:diagram}.
\hspace*{\fill}$\Box$\\
\par
For completeness we include the following well-known proposition.
\begin{Lem}
\label{lem:corners_exp}
All facets of $K$ are exposed faces of $K$, all corners of $K$ are
exposed points of $K$.
\end{Lem}
{\em Proof:\/}
Indirectly, if a face $F$ of $K$ is not exposed, then 
$F\subsetneq{\rm cl}_{\mathcal E}(F)$. Now
${\rm dim}(F)<{\rm dim}\,{\rm cl}_{\mathcal E}(F)<n$ follows by
\cite{Rockafellar} Corollary 18.1.3 and Lemma 4.6 in
\cite{Weis_cones}. This shows that $F$ is not a facet.
\par
Let $F$ be a corner of $K$. First, the face $F$ is exposed: 
Its normal cone is $\Psi(F)=\Psi\circ{\rm cl}_{\mathcal E}(F))$
by Proposition~\ref{pro:diagram}.
By contradiction, if $F\subsetneq {\rm cl}_{\mathcal E}(F)$, then
${\rm cl}_{\mathcal E}(F)$ contains a segment and its normal cone
has dimension $\leq n-1$. This shows that $F$ is an exposed face.
Second, let $x$ belong to the relative
interior of $F$, then $F=\{x\}$: By Definition~\ref{def:varphiN}.2 
of $\Psi$ we have $N(x)=\Psi(F)$ and the proper inclusion
$\{x\}\subsetneq F$ leads to a contradiction as before.
\hspace*{\fill}$\Box$\\
\begin{Cor}
\label{cor:facets_corners}
The conjugate face map
$\mathcal C_K:\mathcal F_K\to\mathcal E_{K^\circ}$ restricts to a 
bijection $\{\text{ facets of }K\;\}\to\{\text{ corners of }K^\circ\;\}$.
The inverse map is the restriction of 
$\mathcal C_{K^\circ}:\mathcal F_{K^\circ}\to\mathcal E_K$ to the
bijection
$\{\text{ corners of }K^\circ\;\}\to\{\text{ facets of }K\;\}$.
\end{Cor}
{\em Proof:\/}
Since facets and corners are exposed faces by Lemma~\ref{lem:corners_exp},
we can use the decomposition of
$\mathcal C_K|_{\mathcal E_K}:\mathcal E_K\to\mathcal E_{K^\circ}$ into the 
bijections
$\mathcal C_K|_{\mathcal E_K}=\Phi\circ{\rm pos}|_{\mathcal E_K}$
in Proposition~\ref{pro:diagram}. Now it suffices to notice from 
(\ref{eq:inv_pos}) that ${\rm pos}|_{\mathcal E_K}$ restricts to a
bijection between the facets of $K$ and the normal cones
($\neq\mathbb R^n$) of $K^\circ$ of dimension $n$.
\hspace*{\fill}$\Box$\\
\par
We arrive at our main results. Let $(L,\leq_L)$ be a poset with greatest 
element $1$. An element $x\in L$, $x\neq 1$ is a {\it coatom} of $L$ if
for all $y\in L$ the two conditions $x\leq_L y$ and $y\neq 1$ imply $y=x$.
We consider for normal cones $N\in\mathcal N_K$ the principal ideal
\[\textstyle
\mathcal T_K(N):=\{T\in\mathcal T_K\mid
T\leq_{\mathcal T}N\}\,.
\]
It is clear that $\mathcal T_K(N)$ is a complete sublattice of the 
touching cone lattice $\mathcal T_K$ and that $N$ is the greatest element 
in $\mathcal T_K(N)$. We call the normal cone $N$ {\it complete} if all 
coatoms of the ideal $\mathcal T_K(N)$ are normal cones of $K$. Otherwise 
$N$ is {\it incomplete}. We also consider for exposed faces
$F\in\mathcal F_K$ the principal ideal
\[\textstyle
\mathcal F_K(F):=\{G\in\mathcal F_K\mid G\leq_{\mathcal F}F\}\,,
\]
which is the face lattice of $F$.
\begin{Thm}
\label{thm:non-exposed_incomplete}
The conjugate face map $\mathcal C_K:\mathcal F_K\to\mathcal E_{K^\circ}$ 
restricts to a surjective map $\mathcal F_K\setminus\mathcal E_K\to
\{F\in\mathcal E_{K^\circ}\mid F\text{ has an incomplete normal cone }\}$. 
It restricts further to a surjective map with range
$\mathcal C_K(\mathcal F_K\setminus\mathcal E_K)$ and with domain equal
to those non-exposed faces of $K$ which are coatoms of $\mathcal F_K(F)$
for some exposed face $F$ of $K$. The preimage of 
$F\in\mathcal E_{K^\circ}$ under $\mathcal C_K$ is
$\mathcal C_K^{-1}(F)={\rm cl}_{\mathcal E}^{-1}(\mathcal C_{K^\circ}(F))$.
\end{Thm}
{\em Proof:\/}
We use Proposition~\ref{pro:diagram} extensively in the proof.
About the preimage of an exposed face $F\in\mathcal E_{K^\circ}$ we 
notice for faces $G\in\mathcal F_K$ that
\[\textstyle
\mathcal C_K(G)\;=\;F
\quad\iff\quad
\mathcal C_{K^\circ}\circ\mathcal C_K(G)\;=\;\mathcal C_{K^\circ}(F)
\quad\iff\quad
{\rm cl}_{\mathcal E}(G)\;=\;\mathcal C_{K^\circ}(F)\,.
\]
\par
We prove that the conjugate face of any non-exposed face has an
incomplete normal cone. For a non-exposed face $F\in\mathcal F_K$ we
consider the touching cone $T:={\rm pos}(F)\in\mathcal T_{K^\circ}$
and we consider the normal cone of its conjugate face
\[\textstyle
N
\;:=\;
\Psi\circ\mathcal C_K(F)
\;=\;
{\rm cl}_{\mathcal N}\circ{\rm pos}(F)
\;=\;
{\rm cl}_{\mathcal N}(T)\,.
\]
Since $F\lneq_{\mathcal F}{\rm cl}_{\mathcal E}(F)$ the lattice
isomorphism $\mathcal F_K\to\mathcal T_{K^\circ}$ of the positive hull 
operator ${\rm pos}$ implies
$T\lneq_{\mathcal T}{\rm cl}_{\mathcal N}(T)=N$. By
{\it Hausdorff's Maximal Principle} there
exists a maximal chain $C$ in the ideal $\mathcal T_{K^\circ}(N)$
including $T$ and $N$, see Chap.\ VIII.7 in \cite{Birkhoff}. A proper
inclusion $F_1\subsetneq F_2$ of faces of $K$ implies a dimension
difference ${\rm dim}(F_1)<{\rm dim}(F_2)$ by Corollary 18.1.3 in 
\cite{Rockafellar}. Hence every chain in the face lattice
$\mathcal F_K$ is finite and hence every chain in the touching cone
lattice $\mathcal T_{K^\circ}$ is finite. So the penultimate element
$P$ in $C$ exists and $P$ is a coatom in $\mathcal T_{K^\circ}(N)$
because $C$ is a maximal chain. By contradiction, if $N$ is a complete 
normal cone, then $P\in\mathcal N_{K^\circ}$. Then 
$T\leq_{\mathcal T}P\lneq_{\mathcal T}N$ implies
${\rm cl}_{\mathcal N}(T)\leq_{\mathcal N}P$ and this contradicts
${\rm cl}_{\mathcal N}(T)=N$. 
\par
We prove surjectivity for the second, smaller, restriction. It suffices to
find for every exposed face $F\in\mathcal E_{K^\circ}$ with incomplete
normal cone $N:=\Psi(F)$ a coatom $G$ of
$\mathcal F_K(\mathcal C_{K^\circ}(F))$ which is a non-exposed face of $K$
and to show $\mathcal C_K(G)=F$. There exists a coatom
$T\in\mathcal T_{K^\circ}(\Psi(F))$ such that 
$T\not\in\mathcal N_{K^\circ}$ and we put $G:={\rm pos}^{-1}(T)$. Since 
${\rm pos}$ is a lattice isomorphism
$\mathcal F_K\to\mathcal T_{K^\circ}$, the face $G$ is a coatom of
$\mathcal F_K(\mathcal C_{K^\circ}(F))$. Since ${\rm pos}$ 
restricts to a bijection ${\rm pos}:\mathcal E_K\to\mathcal N_{K^\circ}$
from the exposed faces to the normal cones, $G$ is a non-exposed face of
$K$. Finally
\[\textstyle
\mathcal C_K(G)
\;=\;
\Phi\circ{\rm pos}(G)
\;=\;
\Phi(T)
\;=\;
\Phi\circ{\rm cl}_{\mathcal N}(T)
\;=\;
\Phi(N)
\;=\;
F   
\]
follows because ${\rm cl}_{\mathcal N}(T)=N$ holds as 
$T$ is a coatom of $\mathcal T_{K^\circ}(N)$.
\hspace*{\fill}$\Box$\\
\begin{Rem}
\label{rem:tothm}
\begin{enumerate}[1.]
\item
A face with an incomplete normal cone is a singular face (with 
normal cone of dimension $\geq 2$). Indeed, a one-dimensional normal
cone is a closed ray $r$ and its two non-empty faces $\{0\}$ and $r$
are both normal cones.
\item
If we apply the dual lattice isomorphism $\mathcal C_{K^\circ}$ to
the second restriction in Theorem~\ref{thm:non-exposed_incomplete} 
then it says that the preimage ${\rm cl}_{\mathcal E}^{-1}(F)$ of an
exposed face $F$ of $K$ under the exposed closure
${\rm cl}_{\mathcal E}$ contains a coatom of the face lattice of $F$
whenever ${\rm cl}_{\mathcal E}^{-1}(F)\supsetneq\{F\}$.
\par
In higher dimensions $n\geq 4$, of course, ${\rm cl}_{\mathcal E}^{-1}(F)$
can contain non-exposed faces, which are not coatoms of the face lattice
of $F$. An example is the direct sum of two copies of the convex body in 
Figure~\ref{fig:dual_pics} a) or c).
\end{enumerate}
\end{Rem}
%
%
%
%
%
\section{Conjugate faces in dimension two}
\label{sec:2D}
\par
We study conjugate faces of a polar pair of planar convex bodies, in 
particular we study the conjugate faces of non-exposed points and we
count special points of the two convex bodies. We characterize self-dual
planar convex bodies without non-exposed faces and we provide a
general construction for them. 
\par
Let $K\subset\mathbb R^2$ be a convex body with $0\in{\rm int}(K)$ and 
polar convex body $K^\circ$. For extremal points $x\in K$ there are two 
alternatives. They have a normal cone $N(x)$ of dimension
${\rm dim}N(x)=1$ resp.\ ${\rm dim}N(x)=2$,
\[\textstyle
x \text{ is smooth}
\qquad\text{resp.}\qquad
x \text{ is a corner.}
\]
The normal cone of a corner $x\in K$ is a {\it salient}\footnote{%
If $N(x)$ contains a line, then $K$ is included in a hyperplane in 
$\mathbb R^2$, see e.g.\ (15)(iv) in \cite{Weis_cones}, and 
${\rm int}(K)=\emptyset$ follows.} convex cone i.e.\ 
a convex cone such that $N(x)\cap(-N(x))=\{0\}$. It follows that
$N(x)$ has two distinct one-dimensional rays $r_1,r_2$ as its faces.
Three types of corners can be distinguished:\\[1ex]
\centerline{%
\begin{tabular}{ll}
$x$ is a {\it polyhedral} corner if & 
$r_1,r_2\in\mathcal N_K$,\\
$x$ is a {\it mixed} corner if & $r_1\in\mathcal N_K$
or $r_2\in\mathcal N_K$ but not both,\\
$x$ is a {\it free} corner if & 
$r_1,r_2\not\in\mathcal N_K$.
\end{tabular}}\\[1ex]
All facets of $K$ are one-dimensional, we call them {\it segments}. If
an extremal point $x\in K$ lies on a segment $s\subset K$ we call $x$
and $s$ {\it incident}. Any non-empty face of $K$ is either an extremal 
point, a segment or $K$ itself. The extremal points and relative
interiors of segments are a partition of the boundary $\partial K$,
see Theorem 18.2 in \cite{Rockafellar}. The boundary $\partial K$ is 
homeomorphic to the unit circle $S^1$ under a positively homogeneous map 
(see the Theorem of Sz.\ Nagy in Construction~\ref{con:support_radial}).
\begin{Rem}[Local classification of extremal points]
\label{rem:local_corner}
\begin{enumerate}[1.]
\item
Segments and corners are exposed faces by Lemma~\ref{lem:corners_exp}
and the proof that every non-exposed point is incident with a unique 
segment is given in Remark 1.1 in \cite{Weis_cones}. This shows
\begin{align*}
\mathcal F_K&\setminus\mathcal E_K
\;=\;
\{\text{ non-exposed points }\}\\
& \;=\;
\{\text{ smooth extremal points incident with a unique segment } \}
\end{align*}
except the inclusion ``$\supset$'' in the second equality. This 
follows by contradiction from the dual lattice isomorphism 
$\Psi|_{\mathcal E_K}:\mathcal E_K\to\mathcal N_K$ between exposed faces
and normal cones in (\ref{eq:l_anti_normalcone}): If $x$ is an 
exposed point incident with a segment $s$, then $x\subsetneq s$ 
shows $N(x)\supsetneq\Psi(s)$. Then ${\rm dim}N(x)=2$ so $x$ is not
smooth.
\item
By the dual lattice isomorphism
$\Psi|_{\mathcal E_K}:\mathcal E_K\to\mathcal N_K$, exposed points 
$x\in\mathcal E_K$ split into the three types of corners above and
into smooth exposed points:\\[3mm]
\centerline{\begin{tabular}{|lll|}
\hline
$x$ is a polyhedral corner & $\iff$ & 
$x$ is the intersection of two segments,\\\hline
$x$ is a mixed corner & $\iff$ & 
$x$ is incident with a unique segment,\\\hline
$x$ is a free corner or & $\iff$ & 
$x$ is not incident with a segment.\\
a smooth exposed point & & \\\hline
\end{tabular}}\\[3mm]
(We have seen in part 1 that a smooth exposed point is not incident with
any segment.) Examples are depicted in Figure~\ref{fig:dual_pics} 
and~\ref{fig:self_pics}. 
\end{enumerate}
\end{Rem}
\par
To understand the conjugate face map we divide the non-exposed points in
\[\textstyle
\mathcal F_K^{\rm mixed}
\;:=\;
\{x\in\mathcal F_K\setminus\mathcal E_K\mid
\mathcal C_K(x)\text{ is a mixed corner of }K^\circ\,\}
\]
and 
\[\textstyle
\mathcal F_K^{\rm free}
\;:=\;
\{x\in\mathcal F_K\setminus\mathcal E_K\mid
\mathcal C_K(x)\text{ is a free corner of }K^\circ\,\}\,.
\]
We show that $\{\mathcal F_K^{\rm mixed},\,\mathcal F_K^{\rm free}\}$
is a partition of the non-exposed points
$\mathcal F_K\setminus\mathcal E_K$.
\begin{Lem}
\label{lem:easy}
The conjugate face map
$\mathcal C_K:\mathcal F_K\to\mathcal E_{K^\circ}$ restricts to the
surjection
\begin{align*}
\mathcal F_K\setminus\mathcal E_K
\;\longrightarrow\;
\{\text{ mixed corners of } K^\circ\,\}
\;\cup\;
\{ \text{ free corners of }K^\circ\,\}\,.
\end{align*}
The restriction of $\mathcal C_K$ to $\mathcal F_K^{\rm mixed}$ is
$1$-$1$ and the restriction to $\mathcal F_K^{\rm free}$ is $2$-$1$.
\end{Lem}
{\em Proof:\/}
All one-dimensional normal cones of $K^\circ$ are complete. This shows
\begin{align*}
\{F\in\mathcal E_{K^\circ}\mid
&F\text{ has an incomplete normal cone }\}\\
&\;=\;
\{\text{ mixed corners in }K^\circ\,\}
\,\cup\,
\{\text{ free corners in }K^\circ\,\}
\end{align*}
so Theorem~\ref{thm:non-exposed_incomplete} proves the first claim.
Proposition~\ref{pro:diagram} and the positive hull isomorphism
${\rm pos}:\mathcal F_K\to\mathcal T_{K^\circ}$ show for mixed and 
free corners $x\in K^\circ$ 
\[\textstyle
{\rm pos}\circ\mathcal C_K^{-1}(x)
\;=\;
{\rm cl}_{\mathcal N}^{-1}(N)
\;=\;
\left\{\begin{array}{ll}
\{N,T\} & \text{if $x$ is a mixed corner,}\\
\{N,T_1,T_2\} & \text{if $x$ is a free corner,}
\end{array}\right.
\]
where $N:=N(x)$ is the normal cone and $T,T_1,T_2$ are rays such that 
$T_1\neq T_2$. The inverse (\ref{eq:inv_pos}) of ${\rm pos}$ gives
\[\textstyle
\mathcal C_K^{-1}(x)
\;=\;
\left\{\begin{array}{ll}
\{N\cap\partial K,T\cap\partial K\} & \text{if $x$ is a mixed corner,}\\
\{N\cap\partial K,T_1\cap\partial K,T_2\cap\partial K\} &
\text{if $x$ is a free corner.}
\end{array}\right.
\]
As ${\rm pos}:\mathcal E_K\to\mathcal N_{K^\circ}$ is a bijection
between exposed faces and normal cones, $N\cap\partial K$ is an
exposed face and $T\cap\partial K$ and 
$T_1\cap\partial K\neq T_2\cap\partial K$ are non-exposed points.
\hspace*{\fill}$\Box$\\
\par
We use a $10$-tuple to label the cardinalities (possibly 
$\infty$) of special points and segments:\\
\centerline{\begin{tabular}{l|l|l}
convex body & $K$ & $K^\circ$\\\hline
non-exposed points & $n$ & $n^\circ$\\
polyhedral corners & $p$ & $p^\circ$\\
mixed corners & $m$ & $m^\circ$\\
free corners & $f$ & $f^\circ$\\
segments & $s$ & $s^\circ$\\\hline
\end{tabular}}\\[3mm]
The following linear equations hold for all planar convex bodies
$K$ with $0\in{\rm int}(K)$. From Lemma~\ref{lem:easy} follow the
equations 
\begin{equation}
\label{eq:enum1}\textstyle
n\;=\;m^\circ+2f^\circ\,,
\qquad
n^\circ\;=\;m+2f\,.
\end{equation}
By Corollary~\ref{cor:facets_corners} we have
\begin{equation}
\label{eq:enum2}\textstyle
s\;=\;p^\circ+m^\circ+f^\circ\,,
\qquad
s^\circ\;=\;p+m+f\,.
\end{equation}
Counting endpoints of segments, we get from 
Remark~\ref{rem:local_corner}
\begin{equation}
\label{eq:enum3}\textstyle
2s\;=\;n+2p+m\,,
\qquad
2s^\circ\;=\;n^\circ+2p^\circ+m^\circ\,.
\end{equation}
These equations span a five-dimensional space of linear functionals and
on the other hand the examples in Figure~\ref{fig:dual_pics} plus the
example of a triangle provide five linearly independent data 
vectors.
\par
If $K$ is self-dual then five cardinalities $(n,p,m,f,s)$ suffice to
count the special points and the above equations reduce to two linear 
equations
\begin{equation}
\label{eq:2}
\boxed{s - p \;=\; n - f \;=\; \tfrac 1{2}( n + m )}
\end{equation}
while three linearly independent data vectors are available from 
Figure~\ref{fig:self_pics} a)--c). 
\par
In the following proposition the necessary condition of an odd number of 
segments is likely to be well-known. For completeness we include a proof.
\begin{Lem}
\label{lem:no_non-exposed}
If $K$ is self-dual and has no non-exposed points ($n=0$), then all 
corners of $K$ are polyhedral ($m=f=0$) and $s=p$. Either $K$ is
strictly convex ($s=0$), a polytope with $s=3,5,7,\ldots$ segments 
or $s=\infty$.
\end{Lem}
{\em Proof:\/}
If $n=0$ then (\ref{eq:2}) implies $f=0$ since $n+m\geq 0$. So $m=0$ 
follows and also $s=p$. A two-dimensional convex body without boundary 
segments ($s=p=0$) is strictly convex i.e.\ all boundary points
are smooth exposed points.
\par
We consider $0<s<\infty$. Any endpoint of a segment is an exposed point
(since $n=0$) hence it is a polyhedral corner (since $m=0$). As the
number of segments $s$ is finite, they are connected in a polygonal
circuit. So $K$ is a polytope, which must have at least three 
segments. We show that $p$ is odd.
\par
Like in Construction~\ref{con:assemble_two_half} b) we rotate the
polytope $K$ about the origin such that
$x_-:=-\s{1\\0}\rho_K(-\s{1\\0})$ maximizes the Euclidean norm on $K$. 
Then $x_-$ is an exposed point of $K$, hence a polyhedral corner of $K$ 
and the normal cone $N(x_-)$ is a two-dimensional salient convex cone.
The segments incident with $x_-$ lie in the ball of radius $|x_-|$
about the origin, so $N(x_-)$ meets $\mathbb R^2\setminus H^-$ and 
$\mathbb R^2\setminus H^+$. By (\ref{eq:conj_norm}) we have
\[\textstyle
N(x_-)
\;=\;
{\rm pos}\circ\mathcal C_K(x_-)\,,
\]
so the conjugate face $\mathcal C_K(x_-)$ is a segment meeting
$\mathbb R^2\setminus H^-$ and $\mathbb R^2\setminus H^+$. If
$r:\mathbb R^2\to\mathbb R^2$ is the reflection $a\mapsto(-a)$, then
for some $y\in\mathbb R^2\setminus H^-$ and
$z\in\mathbb R^2\setminus H^+$
\[\textstyle
r\circ\mathcal C_K(x_-)
\;=\;
[y,z]
\]
is a segment of $K^*=K$. In particular, the boundary point
$\s{1\\0}\rho_K(\s{1\\0})$ of $K$ is not a corner since it lies in
the relative interior of $[y,z]$.
\par
We consider the polygonal chain $C$ in the boundary $\partial K$ from 
$x_-$ to $y$ which lies in $H^+$. Its segments are in bijection to its  
vertices distinct from $x_-$ (by assigning endpoints in the direction
from $x_-$ to $y$). As $x_-$ is a corner, the segments of $C$ are the
segments $s\neq r\circ\mathcal C_K(x_-)$ of $K$ meeting
$K\setminus H^-$. As $y\in K\setminus H^-$ the vertices of $C$ 
distinct from $x_-$ are the corners of $K$ lying in $K\setminus H^-$.
This gives a bijection 
\begin{equation}
\label{eq:bij_c_s1}\textstyle 
\{\text{ corners in $K\setminus H^-$ }\}
\;\cong\;
\{\text{ segments $s\neq r\circ\mathcal C_K(x_-)$ meeting
$K\setminus H^-$ }\}
\end{equation}
between a subset of corners of $K$ and a subset of segments of $K$.
Similarly we have
\begin{equation}
\label{eq:bij_c_s2}\textstyle 
\{\text{ corners in $K\setminus H^+$ }\}
\;\cong\;
\{\text{ segments $s\neq r\circ\mathcal C_K(x_-)$ meeting
$K\setminus H^+$ }\}\,.
\end{equation}
By Corollary~\ref{cor:facets_corners} the map $r\circ\mathcal C_K$ 
restricts to a bijection between the segments and the corners of $K$,
one pair of corresponding faces being 
\begin{equation}
\label{eq:bij_c_s3}\textstyle 
\{x_-,r\circ\mathcal C_K(x_-)\}\,.
\end{equation}
The corners and segments (\ref{eq:bij_c_s1}), (\ref{eq:bij_c_s2})
and (\ref{eq:bij_c_s3}) of $K$ are a complete list. Hence, 
according to Lemma~\ref{lem:ind_q}.2~d) the map $r\circ\mathcal C_K$
is a bijection between (\ref{eq:bij_c_s1}) and (\ref{eq:bij_c_s2}).
\hspace*{\fill}$\Box$\\
\par
We provide a construction for planar self-dual convex bodies without
non-exposed points. Part b) shows that the construction is general. 
\begin{Con}
\label{con:2}
\begin{enumerate}[a)]
\item
Let $K$ satisfy $\rho_K(\pm\s{1\\0})=h_K(\pm\s{1\\0})=e^{\pm\lambda}$
for some $\lambda\in\mathbb R$. Let $K$ have no non-exposed points and
let all corners of $K$ be polyhedral. We assume that
$x_-:=-\s{1\\0}\rho_K(-\s{1\\0})$ is a smooth extremal point of $K$ 
if and only if $x_+:=\s{1\\0}\rho_K(\s{1\\0})$ is a smooth extremal 
point of $K$. Then $L:=(K\cap H^+)\cup(K^*\cap H^-)$ is a self-dual
convex body without non-exposed points.
\item
For every planar self-dual convex body $K$ without non-exposed points
exists a rotation in $\psi\in SO(2)$ such that $\psi(K)$ satisfies the 
assumptions in a).
\end{enumerate}
\end{Con}
{\em Proof:\/}
To prove b) we consider a rotated convex body $K$ according to 
Construction~\ref{con:assemble_two_half} b). As rotation is an isometry,
$K$ is self-dual and has no non-exposed faces. All corners of
$K$ are polyhedral by Lemma~\ref{lem:no_non-exposed}. Since $K$ has
maximal $x$-extension on the $x$-axis we have
$x_+\in r\circ\mathcal C_K(x_-)$. If $x_-$ is a smooth exposed point
then Corollary~\ref{cor:smooth_exp} shows that 
$x_+=r\circ\mathcal C_K(x_-)$ is a smooth exposed point and
{\it vice versa}.
\par
We prove a). Construction~\ref{con:assemble_two_half}~a) already shows
that $L$ is a self-dual convex body. We show that $L$ has no
non-exposed points. First we show that any extremal point $x$ of $L$
in $L\setminus H^-$ is an exposed point of $L$ (the case
$x\in L\setminus H^+$ is analogous). By Lemma~\ref{lem:ind_q}.1 $x$ is
an extremal point of $K$, hence an exposed point of $K$. By 
Lemma~\ref{lem:ind_q}.2 c) there exists $u\in\mathbb R^2\setminus H^-$ 
such that $\{x\}=K\cap H_K(u)$. Then Lemma~\ref{lem:ind_q}.2 b) shows
that $\{x\}=L\cap H_L(u)$ is an exposed point of $L$.
\par
We show that a non-exposed point $x_-$ in $L$ leads to a contradiction,
the proof for $x_+$ is analogous. We will use that all corners of $K^*$ 
are polyhedral and (since $(K^*)^*=K$) that $K^*$ has no non-exposed 
points (this is proved in Lemma~\ref{lem:easy}). By 
Remark~\ref{rem:local_corner} any extremal point of $K$ or $K^*$ is
either a smooth exposed point or a polyhedral corner. 
\par
Since $L$ has maximal $x$-extension on the $x$-axis we have
$x_+\in r\circ\mathcal C_L(x_-)$. If $x_-$ is a 
non-exposed point then Lemma~\ref{lem:easy} shows that
$r\circ\mathcal C_L(x_-)$ is a mixed or a free corner of $L$ so
\begin{equation}
\label{eq:x-+}
x_+
\;=\;
r\circ\mathcal C_L(x_-)\,.
\end{equation}
The contradiction that $x_+$ is incident with two segments
in $L$ completes the proof.
\par
If $x_-$ is a non-exposed point of  $L$ then $x_-$ is incident with a 
segment $[x_-,y]$ of $L$ say for $y\in L\setminus H^-=K\setminus H^-$
(the case $y\in L\setminus H^+$ is analogous by arguing with $K^*$ in 
place of $K$). Since $x_-$ is incident with a unique segment, the
smallest exposed face (\ref{def:clvarphi}) of $L$ containing $x_-$ is
the segment
\begin{equation}
\label{eq:cl_n}
{\rm cl}_{\mathcal E}(x_-)
\;=\;
[x_-,y]\,. 
\end{equation}
\par
We show that $x_+$ is incident with a segment of $L$ in $H^-$. The 
extremal point $y$ of $L$ is an extremal point of $K$ by 
Lemma~\ref{lem:ind_q}.1. Hence $y$ is a polyhedral corner of $K$ and
also of $L$. Corollary~\ref{cor:facets_corners} shows that the face
$s:=r\circ\mathcal C_L(y)$ is a segment of $L$. By (\ref{eq:x-+}), by
the equation
$\mathcal C_L(x_-)=\mathcal C_L\circ{\rm cl}_{\mathcal E}(x_-)$ from
Proposition~\ref{pro:diagram}, by (\ref{eq:cl_n}) and since
$r\circ\mathcal C_L$ is antitone, we obtain that $x_+$ is incident
with the segment $s$,
\[\textstyle
x_+
\;=\;
r\circ\mathcal C_L(x_-)
\;=\;
r\circ\mathcal C_L([x_-,y])
\;\subset\;
s\,.
\]
Lemma~\ref{lem:ind_q}.2~d) shows $s\subset H^-$. 
\par
We find a segment of $L$ in $H^+$ incident with $x_+$. If $x_+$ is a 
smooth exposed point of $K$ then $x_-$ is a smooth exposed point of $K$
by Corollary~\ref{cor:smooth_exp}. This is wrong as $x_-$ is incident
with the segment $[x_-,y]$ of $L$ hence is included in a segment of $K$.
Otherwise if $x_+$ is not a smooth exposed point of $K$ it is a 
polyhedral corner of $K$ or lies in the relative interior of a segment
of $K$. In both cases $x_+$ is included in a segment of $K$ meeting
$K\setminus H^-$, hence is incident with a segment of $L$ included in
$H^+$.
\hspace*{\fill}$\Box$\\
\par
We give an example of a planar self-dual convex body with
$n=0$ and $s=p=\infty$.
\begin{Exa}
\label{exa:infty}
Let $E:=\{u(\alpha)\mid\alpha=(\tfrac k{2}\pm 2^{-m})\pi, 
k\in\{1,3\},m\in\mathbb N\}
\cup\{\alpha(\tfrac \pi{2}),\alpha(\tfrac 3{2}\pi)\}$ for
$u(\alpha):=\s{\cos(\alpha)\\\sin(\alpha)}$. It follows from
Carath\'eodory's theorem, see e.g.\ Theorem 17.2 in \cite{Rockafellar},
that the convex hull $K$ of $E$ is compact. Since $E\subset S^1$, the 
convex body $K$ has no non-exposed points. The two accumulation points 
$\alpha(\tfrac \pi{2})$ and $\alpha(\tfrac 3{2}\pi)$ of
$E$ are approximated by points of $E$ both counterclockwise and clockwise
on $S^1$, hence they are smooth exposed points of $K$. This shows that all 
corners of $K$ are polyhedral. The convex body
$(K\cap H^+)\cup(K^*\cap H^-)$ is self-dual and has no non-exposed points 
by Construction~\ref{con:2} a), it is depicted in 
Figure~\ref{fig:self_pics} d).
\end{Exa}
%
%
%
%
%
%
\bibliographystyle{amsalpha}

\end{document}